\documentclass[12pt]{article}
\usepackage[numbers]{natbib}
\usepackage{amssymb}
\begin{document}

\title{\bf The Golden mean, scale free extension of Real 
number system, fuzzy sets and  $1/f$ spectrum 
in Physics and Biology}

\author{
Dhurjati Prasad Datta\\
Department of Mathematics \\
North Eastern Regional Institute of
Science and Technology\\
Itanagar-791109, Arunachal Pradesh, India} 
\date{}
\maketitle
\begin{center}email:dp${_-}$datta@yahoo.com\end{center}

\begin{abstract}
We show that the generic $1/f$ spectrum problem acquires  a 
natural explanation in a class of scale free solutions to the ordinary differential 
equations. We prove the existence and uniqueness of this class of solutions 
and show how this leads to a nonstandard, fuzzy  extension of the ordinary framework of 
calculus, and hence, that of the classical dynamics and  quantum mechanics. 
The exceptional role of the golden mean irrational number is also explained.
\end{abstract}
\begin{center}PACS No.: 02.30.Hq; 05.45.+j; 47.53.+n\end{center}
\begin{center}Key words: golden mean, scale free, stochastic fluctuations, fuzzy sets, 
1/f spectrum \end{center}
\begin{center} {\em Chaos, Solitons \& Fractals, vol 17, no.4, 781-788, (2003)}
\end{center}

\newpage
\section{Introduction}
\par The ubiquitous presence of $1/f$ spectrum in a wide range of Natural and biological 
processes\cite{pr, sc} is considered to be an interesting problem in theoretical 
physics. 
Any new insight into the problem is expected to enhance our understanding of 
the origin of complex structures in Nature.  
In this paper, we present a body of new mathematical results which is likely to 
initiate a new  approach in understanding the origin of complexity 
and the generic $1/f$ spectrum problem in Nature. 
\par In a recent paper~\cite{dp}, we have initiated an investigation which explores the 
role of time inversion in the context of a linear differential dynamical system. 
By `time inversion' we mean the following. We assume that time has a stochastic 
element in its small scale behaviour, viz.; time can `fluctuate' following the 
inversion rule $t_{-}=1/t_{+}$, where $t_{-}=1-\eta,\,t_{+}=1+\eta,\,0<\eta<<1$, 
and $\eta$ may be a stochastic variable, in the neighbourhood of $t=1$, say. 
The ordinary linearly flowing (reversible) time sense of order $t\sim$O(1), then , 
arises only in the mean, when the small scale stochastic fluctuations are properly 
{\em coarse grained}. Our investigations under this assumption then lead to 
the uncovering of an exact class of scale free solutions to the linear equation   
\begin{equation}\label{so}
{\frac{{\rm d}\ln T}{{\rm d}t}}=1
\end{equation}

\noindent We show that thanks to this exact class of solutions, the real number 
system $R$ is identified with its nonstandard extension~\cite{ns} ${\bf R}$, viz.; 
${\bf R}=R$. Consequently, every real number is identified with an equivalence 
class of scale free fluctuations, undergoing stochastic, irreversible evolutions 
following the cascades of the infinite continued fraction of the golden mean 
$\nu = (\sqrt 5-1)/2$. We show that this scale free fluctuation of the stochastic 
fractal time endows a linear system with the generic low frequency $1/f$ 
spectrum. Further, the generic probability distribution of these scale free 
fluctuations turns out to be a Gamma distribution. 

\par The relevance of the present work may be judged in the light 
of the recent works~\cite{mw,pl} uncovering new relationships 
between time and the number theory. 
In ref\cite{sl}, Selvam, on the other hand, developed a cell dynamical model  
for the computational error growth, analogous to  the formation 
of turbulent eddies.  It is shown 
that the ratio of higher to lower precision error domains scales, in a long 
computer run, as the golden mean $\nu$, revealing a universal self similar 
pattern in the round- off error structure. Further, the theory of Cantorian spacetime of 
El Naschie [8-12] seems to anticipate some of the features of the scale free fractal 
time.
 
\par In the present paper, we reverse our previous approach and base our analysis 
directly on the exact scale free solutions of eq(1). Our main results in the form of 
certain 
theorems are proved in the next section. We show that the real number system is a 
fuzzy, nonstandard set, infinitesimally small elements of which experience a stochastic, 
irreversible evolution. 
In Sec.3, we show how these intrinsic fluctuations  
lead to the generic 1/$f$ spectrum for a general differential dynamical system. We 
close in Sec.4 with some concluding remarks.

\section{Main theorems and their proofs}

Our main results are stated in the following three theorems.
\newtheorem{theo}{Theorem}
\begin{theo}
There exists a unique class of one parameter family of scale free 
solutions 
\begin{equation}\label{ns}
\ln T(t)=t + k {\frac{T(\ln t)}{t}},\, t\neq 0
\end{equation}

\noindent $k$ being an arbitrary parameter, to eq(1). These solutions 
are fundamentally non-local. 
\end{theo}

Consequently, 
\begin{theo}
This non-locality directs 
a real variable $t$ to change, in the neighbourhood of $t=1$,  by a `local inversion', 
$t_{-}=1/t_{+}$, where $t_{-}=1-\eta,\,t_{+}=1+\eta,\,0<\eta<<1$, 
in contrast to the conventional mode of linear increments. Consequently, the solutions 
are only continuously first order differentiable, at the points where $t$ changes by 
local 
inversions.
\end{theo}
   
\begin{theo}
 The family of solutions is also shown to provide a `dynamical' representation 
of the real number system, each member of which 
is identified with an equivalence class of scale free fluctuations, undergoing 
a stochastic, irreversible evolution, following the cascades of of infinite 
continued fraction of the golden mean. 
Consequently, the real number set $R$ 
is a fuzzy, nonstandard set.

\end{theo}

\par The existence part of Theorem 1 is almost trivial. To verify that 
eq(2) is indeed a solution, we consider instead a more general ansatz 
$\ln T(t)=t + k{\frac{\tau(t)}{t}}$, and assume that both $T$ and $\tau$ 
are continuously first order differentiable in $t$. A direct differentiation now 
yields ${\frac{{\rm d}\ln T}{{\rm d}t}}-1={\frac{k}{t^2}}(-\tau+t{\frac{{\rm d}\tau}
{{\rm d}t}})$, which is true for all $t\,(\neq 0)$ (and $k$). One thus obtains 
eq(1), along with the fact that $\tau(t)\equiv T(\ln t) $, since 
\begin{equation}\label{ss}
t{\frac{{\rm d}\ln \tau}{{\rm d}t}}=1
\end{equation}

\noindent  We note that $T(t)$ is indeed a new solution, 
since it is undefined at $t=0$. Clearly, the standard solution of 
eq(\ref{so}) $\ln T_s=t,\,T_s(0)=1$  follows from 
the above ansatz provided we choose $k=0$. However, in the present 
framework, there is no reason for this `a priori' restriction.  
The form of $\tau$ means that when $T$ is a function, for instance,
 of $t\gtrapprox 1$, then $\tau(t)$ must be a self similar replica of  $T$ 
on the smaller logarithmic scale $\ln t$, so that $\tau(t\approx 1)=T(t\approx 0)$. 
The total solution $T(t)$ is thus obtained as a co-operative effect of two different 
scales that arise naturally in eq(1), viz., when the variable $t$ is allowed to 
change from $t\gtrapprox 0$ to $t\approx 1$, and then is mapped back to 
 $t\gtrapprox 0$ on the logarithmic scale: $dt\rightarrow d\ln t$ near 
$t\approx 1$ (c.f., remark 1). Indeed, eq(\ref{ss}) coincides with eq(1) in the limit 
$t\rightarrow 1$.
If one wishes, instead,  to write $\tau=ct,\,t>0$ as usual, then eqs(1) and (3) 
match over the {\em ordinary scale} $t\sim$O(1), $\ln t\sim$O(1), 
viz., $1\lessapprox t \sim $O($e$), provided $k\approx 0,\, c\approx 1$, 
when the initial condition is defined at $t=1:\, \ln T(1)=1$ ( imposing initial 
condition at $t=0$ is not allowed since $t\neq 0$). The standard solution 
$T_s$ is thus realized as a particular solution, when this more general class 
of solutions is restricted to the ordinary scale $t\sim 1$. In anticipation 
of the following, we call the nontrivial part ($T_f$) in $T$ as the `fluctuation' over 
the standard mean solution $T_s$. Further, the solution is, in general,  non-local, 
because $T(t)$ defines a non-local (`distant') connection between two logarithmically 
separated points on the $t$- axis, for instance, $t=1$ and $t=0$ 
(equivalently,$t\rightarrow 0^+$ and $t\rightarrow \infty$).  The non-locality drops out 
(approximately ) only over the ordinary scale. 
\newtheorem{rem}{Remark}

\begin{rem}
{\rm We note that  
the standard general solution of eq(1), but restricted to the punctured real line 
$R-\{0\}$, viz.; $\ln T=t + c, \, c$ a constant number, and $t\neq 0$, 
could be interpreted  to belong to the generalised class of solutions of eq(1),  
for an unrestricted $t$, and vice versa, since 
$\ln T(t)=t + k{\frac{\tau(t)}{t}},\, c=kc_0$ and $\tau=c_0t$ is the standard solution 
of 
eq(\ref{ss}).  Because of the translation symmetry of eq(1), the initial value 
problem in either of the cases is defined only approximately 
at the initial point $t\approx 0$. The constant $c$ is thus determined only 
approximately, which then acts as a seed for the residual self similar evolution 
for the supposedly constant fluctuating component $T_f$.}
\end{rem}
  
\begin{rem}
{\rm 
An alternative way of seeing the possibility of a nontrivial fluctuation for a solution 
of eq(1) is to use the Born-Oppenheimer- like factorisation~\cite{dp1}
$T=T_sT_f$ so that eq(1) reduces to 
 $t^{-1}{\frac{{\rm d}T_f}{{\rm d}\ln t}}=0$.  Clearly, for any moderately large value 
of $t(>0)$ $T_f$ is a constant, but for an arbitrarily large $t:\,t^{-1}=\epsilon\approx 
0$, ${\frac{{\rm d}T_f}{{\rm d}\ln t}}$ can indeed be O(1).} 
\end{rem}

\par To study the salient features of the solution, we note that the standard 
solution $\ln T_s=t$ defines a 1-1 
(identity) mapping $R\rightarrow R$. Let $\phi(t)=\tau(t)/t$, $\tau(t)$ being a 
nontrivial solution of eq(\ref{ss}). Then the new solution eq(\ref{ns}) gives an 
extension 
of $T_s$ to a class of scale free solutions, since $\tau$, being a solution 
of eq(\ref{ss}),  is scale free. We note that the fluctuation $\phi(t)$ is also scale 
free, since $\phi$ satisfies the scale free equation

\begin{equation}\label{cns}
t{\frac{{\rm d}\phi}{{\rm d}t}}=0
\end{equation}

\noindent and hence represents a `slowly varying constant'.  
Both $\tau$ and $\phi$ satisfy the scaling law 
\begin{equation}\label{i}
f(kt)=kf(t)
\end{equation}

\noindent where $f(t)$ is a scale free function. The scaling law 
for $\tau$  follows from eq(\ref{ss}). Consequently, 
$k\tau(t)/t=\tau_1(t_1)/t_1,\,t_1=kt, 
\tau_1=k\tau$, and hence $\phi$ also satisfies eq(\ref{i}). It also follows 
that $\phi$ is {\em a `universal' function, defined 
in the vicinity of $t=1$}. Because of 
the translation symmetry $t\rightarrow  t-t_0$ of eq(1),  the fluctuation of 
$T(t)$ is actually defined pointwise $\ln T(t)=
t + k \phi(t^\prime),\, t^\prime=t-t_0$. By 
the above universality, however, $\phi(t^\prime)=\phi(t)$,  and hence 
the `nontrivial' neighbourhood, denoted `$t$', of every real number $t$ is 
identical to `1': `$t$'-$\{t\}$=`1'-$\{1\}$. Here, $\{t\}$ is the singleton 
set of the real number $t$.  Consequently, $\phi(t)$ represents 
the universal scaling behaviour of the neighbourhood of 
every real number.
\par Let $c$ denote a constant solution of 
eq(\ref{cns}), in the ordinary sense. Then a scale free extension is given by 
$\phi_c (t)=c(1+k\phi_c(\ln t))$. Rescaling $\phi\rightarrow \phi_c/c$ and 
$k\rightarrow kc$, and using the scaling law eq(\ref{i}), 
we get  $\phi\in$`1' 
in the form $\phi(t_1)=1+k\phi(\ln t_1)$. Here, 
$t_1$ is another `universal 
local parameter', defined in `1'. In-fact, by the inversion 
symmetry of eq(\ref{ss}), viz.; ${\frac{{\rm d}\ln \tau}{{\rm d}\ln t}}=1
\Leftrightarrow {\frac{{\rm d}\ln t}{{\rm d}\ln \tau}}=1$, 
$t$ can be treated as a function of $\tau$, and thus 
$t_1(t)=t(\tau(t))\in $`1'. The logarithm in the $k$- dependent 
term,  as stated above, signifies non-locality, 
if the lhs is defined at $t_1=1$, then the rhs is defined 
at $\ln t_1=\epsilon\approx 0$ 
and hence at $t_1\approx 1+\epsilon$, by universality.
\par Now, in relation to the scale defined by eq(1), any real number $t$ is 
written in the scale free representation as $t_f(k)=t\phi(t_1)$. Clearly, $t_f$ 
satisfies eq(\ref{ss}), and defines a one parameter family of extensions of the identity 
mapping,  
$1_f: R\rightarrow R$. However,  $1_f=\phi(t_1)$, and hence the scaling parameter $k$ 
must be sufficiently small $0<|k|<<1$. {\em The extension is also unique}, for 
each choice of a fixed, but arbitrarily small,  $k$. For, 
let $\bar \phi(t)$ be another solution of eq(\ref{ss}). Then, $|\phi(t)-\bar \phi(t)|=
|k||\phi(\ln t)-\bar \phi(\ln t)|=\epsilon|\phi(\bar t)-\bar \phi(\bar t)|,\,\bar t
=(k\ln t)/\epsilon$, by the scaling law  eq(\ref{i}), for any small $\epsilon>0$, so 
that 
$\bar \phi=\phi$, again by universality. 
This completes the proof of  Theorem 1.

\begin{rem}
{\rm The nontrivial `slowly varying constant' $\phi(t)$ in the generalised solution 
eq(2) 
tells the impossibility of erasing a {\em residual} $t$ dependence not only 
from the fluctuating component of the solutions of eq(1), but also for any real number. 
In the following, we show that this residual $t$-dependence has an intrinsic time-like 
feature.
In ref\cite{dp}, we have shown 
how {\em the golden mean partition of unity :$\nu^2+\nu=1,\,\nu>0$} could be utilised 
to magnify a practically insignificant time dependence to an observable effect 
over a sufficiently long time scale. 
In fact, the residual fluctuation in log scale, which remains {\em unnoticed } in the 
usual treatment 
of Calculus, gets manifested under a nontrivial SL(2,R) realization of a linear 
translation

\begin{equation}\label{gold}
t=1+ \eta\approx {\frac{1+ \nu^2 \eta}{1-\nu \eta}},
\end{equation}
\noindent for a sufficiently small $0<\eta<<1$. A fluctuation which is negligible 
in the scale of $t\approx 1$ can grow to O(1) when $\ln t\approx \eta\sim$O(1). 
The presence of a residual time dependence 
could also be the reason for the origin of $1/f$ noise in the prime number distribution 
\cite{mw}}.
\end{rem}

\begin{rem}
{\rm In the usual treatment, the scaling law eq(\ref{i}) means $f(t)=tf(1)$. Clearly, 
this is 
derived under the assumption that the equality $kt=1$ is exactly realized. In the 
context 
of the generalised solution eq(2), the validity of an exact equality (as in the standard 
classical analysis) is violated. A key reason of this violation is the non-locality of 
the solution, so that the exact value of $T(t)$ at a particular $t$ is non-computable 
since this requires a string of recursions spanning over a  number of distinct  
logarithmic scales (for instance, $\ln t,\,\ln\ln t,\ldots$ for a sufficiently large 
$t$).
Consequently, the set defined by $T(t)$ is a fuzzy set and hence the real set $R$ 
itself is fuzzy. This proves partially Theorem 4. (For a more precise sense of 
fuzziness, 
see below.)}
\end{rem}

\par Now, to prove Theorem 2, we note that 
$1_f\in $`1', by definition, so that $\phi(t_1)=1_f=t_1$, 
by the uniqueness of the extension.  We thus get 
\begin{equation}\label{ii}
t_1(t)=1+kt_1(\ln t)
\end{equation}

\noindent This also 
justifies the term `slowly varying constant', since for an arbitrarily small $|k|,\, 
t_1\approx 1$, is almost constant over the scale of $t\sim 1$, but undergoes 
(infinitely) slow variation, (c.f., Remark 3) which gets manifest  in the smaller 
scale $\ln t\sim 0$: $dt_1/dt\approx 0$, but $dt_1/d\ln t\sim $O(1).  

\par Now,  to see how {\em a nontrivial definition of 
an inversion follows from eq(2) directly, leading to the 2nd derivative discontinuity,}
let $t_{-}=1-\eta,\,t_{+}=1+\eta$,  for $0<\eta<<1$. Then we have
\begin{equation}\label{iii}
t_1(t_{-})=1+kt_1(\ln t_{-})=1-kt_1(\ln t_{+})
\end{equation}
 
\noindent  
by the scaling law eq(\ref{i}), and 

\begin{equation}\label{iv}
t_1(t_{+})=1+kt_1(\ln t_{+})
\end{equation}

\noindent Hence, even as $t$ 
is assumed to change from $t_{-}$ to $t_{+}$ by translation (say), 
$t_1(t_{-})$ gets linked {\em instantaneously} 
to $t_1(t_{+})$, as in eq(\ref{iii}), because of the non-local connection. 
This instantaneous `distant connection' between two  
distinct points in the $t$ scale tells that the change in the zeroth order scale $t$, 
near $t=1$,  to avoid any paradoxical situation (viz., reaching the 
destination before being started),  must have been accomplished 
instantaneously, instead, by an application of a local 
inversion for a small enough $\eta$. This zeroth order inversion then induces  
analogous mode of inversions even for the small scale variable $t_1$,  
since $t_1(t_{+})=1/t_1(t_{-})$ for $(0<)k<<1$. It thus follows  
that {\em a local inversion must replace the ordinary translation in accomplishing 
changes in the vicinity of a point on the real axis,} for the sake of consistency.
It is now easy to verify the 
2nd derivative discontinuity of the scale free solution of eq(3), by differentiating 
eq(\ref{iii}) and eq(\ref{iv}), and then taking the limit $\eta\rightarrow 0^+$. 
Clearly, the first derivatives matches continuously at $t=1$, but 
there is a mismatch in sign in the second 
derivatives. This proves Theorem 2.

\par We note that this discontinuity is not in contradiction with the Picard's 
theorem\cite{pic}. The framework of {\em the Picard's theorem 
(and, as a matter of fact, that of 
the (real) analysis) gets extended under the scale free properties of real numbers 
and local inversions}. We give a thorough independent analysis of this theorem
elsewhere. 

\par We now proceed to prove Theorem 3. We do this in several steps. First,  
we show why {\em the above  discontinuity ought to be be removed to a 
very small ( large ) value of $t$.} We begin by showing, in an alternative way, 
 how the extension of the real system $R$ defined by the generalized 
mapping $1_f$  {\em forbids} an {\em exact} evaluation of a real 
number (c.f.Remark 3). We note from eq(\ref{ss}) that  
\begin{equation}\label{v}
\tau(t_{-})=\tau(t_{+}^{-1})=1/\tau(t_{+})
\end{equation}

\noindent  which follows by 
replacing $t$ by $t^{-1}$ in the equation and identifying the `-' sign in the two   
alternative ways. Further, 
\begin{equation}\label{vi}
\tau(t_{+})/t_{+}=t_{-}\tau(t_{-}^{-1})= t_{-}/\tau(t_{-})
\end{equation}

\noindent and hence 
\begin{equation}\label{vii}
\phi(t^{-1})=1/\phi(t),\,\Rightarrow \phi(t_{-})=1/\phi(t_{+})
\end{equation}

\noindent Clearly, these relations are valid 
even in $\ln t$. Consequently, we have $t_1(t_{-})=1+kt_1(\ln t_{+}^{-1})=
1+k/t_1(\ln t_{+})$, where $t_1(\ln t_{+})=1+kt_1(\ln(\ln t_{+})^{-1})$. We note that 
near $t=1, \,0<\eta<<1$, so that $t_1(\eta)=1+kt_1(\ln(\eta)^{-1})$ and hence in the 
limit $\eta\rightarrow 0^+$, we get $t_1(0)=1+t_1(\infty)$, by the scaling law 
eq(\ref{i}), 
since $k\infty=\infty$, for {\em any} finite $k$. Further, by definition, 0 and 
$\infty$ are related by an  inversion, so that  $t_1(0)=1+1/t_1(0)$ and 
$1/t_1(\infty)=1+t_1(\infty)$. We thus have $t_1(0)=1+\nu$ and 
$t_1(\infty)=\nu$, where $\nu=(\sqrt 5-1)/2$, is the golden mean irrational 
number. Apparently it means that the {\em exact } value of the golden mean is 
realized as the limit of the universal scale free solution of eq(\ref{ss}), 
as $t$ approaches either $\infty$ or 0, when the inversion  
$t_1(\infty)t_1(0)=1$ is exactly satisfied. In reality, however, {\em the 
realization of the exact value of $\nu$  is impossible, in principle.} 
( We note that the local inversion is valid upto O($\eta^2$).) 

\par Now, the above assertion follows if one assumes 
the existence of a nontrivial {\em infinitesimally 
small number}, in the sense of the nonstandard analysis\cite{ns}. In fact, 
{\em we show that such a number exists}. Let the solution of eq(\ref{ss}) 
be written as $t_f=t+k\phi(t)=t(1+k\phi(t)/t)$.\footnote{Note that the 
scale free extension of R is two-fold, viz.; $t_{f+}=t\phi(t),\,t_{f-}=t\phi(t^{-1})$, 
for each $t\in R$ so that $t_{f+}t_{f-}=t^2$, for each $k>0$.} It follows from 
eqs(\ref{ii},\ref{vii}) that 
\begin{equation}\label{viii}
\phi(t^{-1})=1+\tilde k\phi(t),\,\tilde k=k/t
\end{equation}

\noindent Now, for each 
$0<t<1$, there exists a $k$ so that $0<\tilde k<t$. So, in the sense of a limit 
in the ordinary calculus, $\tilde k\rightarrow 0$, as $t\rightarrow 0$, 
leading to a contradictory result that $(\phi(\infty)=)\nu=1$. This contradiction, 
however, is removed in the context of a `physical' limit, which allows a real 
variable $t$ (say) to approach 0 (say, but {\em never coinciding} with 0), 
by {\em exploring smaller and smaller scales}. 
It turns out that this exploration of scales is a never ending process, and so, 
in principle, {\em delays the realization of the final limit to an infinitely distant 
epoch}. 
We note that limit (and for that matter, any continuous process of change ) is 
basically a {\em dynamical} concept, and must have an intrinsic correlation with 
a sense of time. The continuous (monotonic) change of $t$ to smaller and 
smaller values (scales) can indeed be correlated to an increasing sense of 
time through an iterative process.   
( In ordinary treatment, $t$ traverses 
the interval $[0,1]\downarrow $ in unit time. In the context of the scale free 
solution of eq(1), the relevant interval is (0,1] and as argued below, the 
approach $t\rightarrow 0^+$ should, in principle, be a never ending process  
(c.f., the concluding remarks).) We now show that 
(A) {\em the small scale variable $t_1=(\phi(t))$ is intrinsically stochastic 
with an irreversible arrow of `time'.}

\par Indeed, we note that (1) $\tilde k$ can be  a `slowly varying 
constant', for an {\em infinitesimally small}, slowly varying  $k(t)$ satisfying eq(2). 
Consequently, as $t\rightarrow k$ ($t$ varies faster in comparison to $k$, and 
must cross $k$ in a `finite' elapse of time), the applicability of the ordinary 
sandwich theorem is violated, because of the reversal of the above inequality:  
$(k\sim) t<\tilde k (\sim 1)$, allowing a local inversion to materialise, in the 
neighbourhood of the slowly varying $k$. Further, (2) the scale free representation 
\begin{equation}\label{ix}
\nu=\nu_f/\phi(t^{-1})
\end{equation}

\noindent tells that {\em the golden mean number 
$\nu$ must be identified with a uncountable set of fluctuating, 
approximate values (evaluations) of the said number}. 
Consequently, as $t$ approaches smaller and smaller scales, an approximate 
$t$-dependent value of $\nu$, denoted $\nu_{[t^{-1}]},\,[t^{-1}]$ being the 
greatest integer value function, and given by eq(\ref{ix}), 
now approaches slowly to more and more accurate evaluations of $\nu$ 
through the sequence of convergents $\{\nu_n\}$ of the golden mean continued fraction.    
 Indeed, we have $\nu_n=\nu_f/(1+(k(t)/x_nt)),\,x_n=\nu_f/\nu_n,\, 0<k<<t^2$, 
for a sufficiently large $n=[t^{-1}]$, but  $\nu_n=\nu_f/(1+1/x_n)$, at $t\sim k$, 
so that $\nu_{n+1}=\nu_f/(1+1/(1+(k(t_1)/x_{n+s}t_1)))$, 
where $t_1=kt, \,t$ being a decreasing O(1) variable from $t=1$, and so on. We note that 
the self-similarity of intervals $(k^2,k)$ and $(k,1)$ tells that there exists a 
$0<k_1(t)=
k(t_1)<<t_1^2$ so that $x_n=1+(k(t_1)/x_{n+s}t_1),\, n+s=[t_1^{-1}]\approx [k^{-1}],\,
\nu_{[t_1^{-1}]}=\nu_{n+s}$. It thus follows that, for $\nu_f=\nu_{n-1}$, 
eq(\ref{ix}) represents an intrinsic evolution of the convergents $\nu_n$ from a lower 
precision to higher and higher precision values, as $t$ tends to 0 slower and 
slower, since $dt_1/dt<<1$. Further, any approximate value $\nu_f$ of $\nu$ 
can be rescaled to an convergent: $c\nu_f=\nu_{n-1}$ so that 
$\nu_n=\nu_{n-1}/\phi(ct^{-1})$.

\par Let us note here that $k/t$ remains arbitrarily small, but non-zero, 
even in the limit $n\rightarrow \infty$. Consequently, even as all the ordinary 
positive real numbers, represented by the scale $t$,  are exhausted in the limit 
$t\rightarrow 0^+$, $k/t$ remains dynamically active, although $k=0$ in the limit. 
We thus {\em define an infinitesimally small $k(>0)$} by the condition that 
$[t/k]=s$ be arbitrarily large as $t\rightarrow 0^+$. Clearly, an infinitesimal $k$ 
{\em exists}, by the above construction. Obviously, $K$ is {\em infinitely 
large} if $1/K$ is infinitesimal. Further, infinitely large and 
infinitesimally small numbers are (relative) scale dependent, for instance, 
$k$ is infinitesimal relative to the scale $t\sim $O(1) $\Leftrightarrow k$ is 
infinitely large relative to $t\sim $O$(k^2)$. 
\par Now, returning to the golden mean convergents, we note that 
as $t$ explore the $n$th 
order smaller scale, the unfolding of the $n$th level cascade of the continued 
fraction is activated. Thus $\nu_n$ tends to the exact value $\nu$, as $t^{-1}$ 
approaches $\infty$ exploring {\em infinitely} larger and larger scales. But, 
{\em this unfolding of higher order cascades must be a never ending process}, 
since, by the self-similarity of real axis over scales, there can not be a largest 
infinitely large real number. Thus, {\em any} approximate value of $\nu$ must 
experience an intrinsic evolution towards the exact value, but would never 
be able to attain the final accuracy.  We remark that although there is no exact 
value of $\nu$ in the present formalism, a progressive evolution towards better 
and better accuracy values still makes sense. As stated before,  
any iterative process generates a {\em directed sense of time}. 
However, the exact injection moments of local inversions, leading the evolution 
to multiple scales, is uncertain, because of (a) the arbitrariness in the sign 
of the infinitesimal scaling constant $k$ and (b) the absence of an exact value.  
Hence, the intrinsic time sense is {\em stochastic} as well. This proves 
our assertion (A).\footnote{ The directed evolution is global in the sense that 
every real number evolves from a lower to higher precision values following 
the sequence of the golden mean convergents. This global irreversibility 
should also be addressable in the sense of the informational complexity (entropy)
( or in terms of the Hausdorff dimension in $\mathcal{E}^(\infty)$ spaces 
(see e.g.,M. S. El Naschie, Chaos, Solitons \& Fractals, 2002, 14(7), 1121)). 
We defer this approach for a separate investigation.}
 In the following, 
we, however, disregard the randomness from (b) ( which will be  
considered separately).
\par Further, every real number $r$ can be written as $r=r_0\nu$, for a suitable 
$r_0$,  so that the above analysis applies to every $r\in R$.
 Consequently, the scale free  eq(2) can indeed endow every real number  $r$
with an intrinsically  time dependent, scale free representation $r_f(t)$, which is 
essentially a (continuous) distribution of approximate values (fluctuations)
undergoing a slow, stochastically scale free, cascaded evolution down the 
infinite staircase of the golden mean continued fraction. We note also that 
a real number $r$ belongs to the equivalence class of any other number, for 
instance, 1 say, defined by the scale dependent (interval-valued fuzzy) 
membership function\cite{fz}: $r_f=r\phi(t_1^{\pm 1},k)$. Clearly, an approximate 
evaluation $r_f$ of $r$ belongs to the fuzzy number `$r_f$' with a membership 
value $\phi(t_1^{\pm 1},k)\lesssim 1$, for $t_1\sim 1$ but belongs to the fuzzy `1' with 
a membership $0<\phi(rt_1,k)<1$, 
if $0<r<1$, and $0<1/\phi(rt_1,k)<1$, if $r>1$, for $t_1\sim 1/r$. 
A detailed study of the fuzzy aspects 
will be considered separately. This completes the proof of    
Theorem 3. We note that the nonstandard extension of the real number system 
constitutes a valid model of  analysis\cite{ns}. {\em The self-consistent 
derivation of a nontrivial solution on the basis of the scale free eq(3) therefore 
reveals 
hitherto unexplored new features of the real  number system, thereby elevating the 
possible solutions to eq(1) to the class of finitely differentiable solutions.} We note 
that infinite differentiability on the scale $t$ reappears in the present context, when 
the realisation of inversion is removed to an infinitely distant moment $t=\omega$. 
However, for any $t$ there exists $k$ so that $kt\approx 1$. Thus, this 
postponement of inversion, and consequent scale changes, could not be maintained 
for an indefinitely long period of time, involving many infinitely large scales. 
This proves the {\em inevitability} of the scale free 
solutions in the context of eq(1). 

\section{Power spectrum}

\par We note that the randomness in sign raises $t_f$ to a
 random variable ${\bf t_f}$. Since both the signs are equally likely, 
the expectation value of ${\bf t_f}$ 
is written as $<{\bf t_f}>=t(1+\tilde k(t_1)\phi(t_1)),\,\tilde 
k(t_1)=(k(t_1)+k(1/t_1))/2$, 
which has the same form of a nonrandom $t_f$ with $k(>0)$ replaced by $\tilde k$. 
We note that the small (late) $t$ asymptotic form of the fluctuation in the 
scale free solution eq(2), and that of $<{\bf t_f}>$ is given by  
$\bar t_f\sim t^{ \mu},\,\mu=k{\phi\over \ln t}$, according as 
$t\rightarrow 0$ or $\infty$. Here, $\bar t_f=t_f/t$ denote the (total) scale free 
fluctuation of eq(2), over the scale $t$, when the standard mean 
component is removed. The above universal fluctuation 
pattern is obtained simply by rewriting the ansatz for the scale free 
solution for eq(2).  The large $t$ asymptotic form of 
the correlation function of the fluctuation spectrum  
 is now written as $C(t)=<{\bf \bar t_f\bar t_f(0)}>\approx t_f^2 \sim t^{2\mu}$, 
since $\bar t_f(t)=\bar t_f(1),\,\forall t$, by the universality of the scale 
free solution.
\par We remark that in the framework of the classical measurement 
hypothesis, any real physical variable is exactly measurable. The above 
analysis shows that this hypothesis is violated because of the irreducible  
scale free fluctuations. Further, in a realistic physical application, any physical 
variable $T$ (say) is measurable with a finite degree of accuracy, say $r<<1$, 
only. The infinitesimal scale free fluctuations in the exact equation (3) thus 
correspond to the uncertainties in the measured values of this physical variable. 
In that case, the actual infinitesimally small scaling parameter $k$, behaves as 
a `relative  infinitesimal', viz., when $k$ is any real number smaller than the accuracy  
limit, $k<r$. The scale free dynamical representation of real numbers now tells that 
any relatively negligible fluctuation would grow in time, and would have nontrivial 
dynamical effects, as reflected in the non-zero late time asymptotic exponent $\mu$.  
The power spectrum of this power law asymptotic correlation function has the 
generic form $1/f^{1+2\mu}$. We note that a small, but nonzero, $\mu$ is sufficient 
to generate the generic $1/f$ spectrum. 
\par Finally, a general (differential ) dynamical system, represented by an  
ODE of the form $f(t,x,\dot x,\ddot x,\ldots)=0$, where $t$ is the ordinary 
time $t$, is, in view of eq(3), an equation, essentially, in $t_f$. The evolution of 
the system thus inherits the scale free properties of time 
itself, apart from any special features 
depending on the explicit nonlinearities in the problem, the signature of which should 
be revealed in the form $\mu=\mu_f+\mu_d$, where $\mu_f=k\phi$ is the universal 
component from the scale free time and $\mu_d$ is a model dependent term.  
Thus the universal presence of $1/f$ spectrum in the Natural processes over a long 
time scale is adequately explained in the framework of the scale free solutions of
the linear ODE eq(1). 
\section{Closing remarks}
\par We conclude with the following remarks. 

1. The conventional linear treatment of the limit  
$t\rightarrow 0$(say) overlooks the presence of nonlinear slowly varying scales 
of the form $0<k(t)t<<t<<1$, because of the tacit assumption that $t$ in 
$0<t<<1$ approaches 0 with the uniform speed 1. Consequently, 0 in the 
ordinary real number system corresponds to $0_f=0^+\cup 0^-,\, 
0^\pm=\{\pm k(t)t, 0<<t<<1\}$, 
since the infinitesimal scales are indeed {\em unobservable} in the 
framework of the linear Calculus.   

2. Let $I(t)=(-t,t),\,t\rightarrow 0^+$. Then the cardinality of $I(t)$ is $c$, the 
cardinality 
of continuum. However, cardinality of $I(0)=0$. This abrupt 
discontinuous transition from $c$ to 0 
is not explained in the linear Calculus. The scale free extension $0_f$ of the 
ordinary 0, on the other hand, is a Cantor -like set with scaling exponent 
$\mu$, which, for $t\rightarrow 0^+$ is identified as the uncertainty 
exponent \cite{ot} for the fattened real number set $R_f$. However, 
the scale free extension, being defined by an exact solution of eq(1) 
means that $R_f\equiv R$, as remarked already. 
Thus every point of the real line is structured  
and accommodates both the point-like and continuum properties 
analogous to a Cantor set. The ordinary point -like (non-fuzzy) 
structure arises under an
approximation, viz., when the continuum of nonlinear infinitesimal scales 
are neglected. The small, but non-zero uncertainty exponent gives 
a measure of the limit of accurate evaluation (measurement) of a 
physical quantity; since the uncertainty exponent $\mu$ tells that 
points separated by a distance less than $\mu$ would be indistinguishable.  

3. In a dynamical system, the ordinary time $t$ is obtained 
as the (zeroth order)  ``mean field'' realization of the scale free time $t_f$, 
as the expectation value of the 
infinitesimal scale free fluctuation $k\phi(t)$ is negligible. The ordinary 
evolution of the system would, therefore, 
continue till $t\sim 1/k,\, k$ being the supremum 
of infinitesimal scales. During the transition period (i.e., when $t\rightarrow 
(1/k)^-)$, the system would experience random fluctuations as the scale free 
component (fluctuation) grows to O(1): $\phi(kt)=\phi(t_{1-})\sim 1$. The 
system, however, returns to the ordinary regular pattern of evolution once 
$t$ crosses $1/k$ from left to right  as $\phi(t_1^{-1})=1+k\phi (t_1)$ 
and $t_1 $ tends slowly (compared to $t$) to $(1/k)^-$. The scale free evolution 
of a system would thus resemble, generically, to an intermittent process as 
$t$ tends to $\infty$ exploring longer and longer scales. Any dynamical 
system, when evolved over an `infinitely' long period of time, would 
inherit naturally this universal pattern of fluctuations of the scale free time.  
Clearly, this universal fluctuation is obtained as a consequence of a 
cooperative effect of random infinitesimal scales, 
which contribute gradually as the evolution is continued over longer 
and longer periods of time. Consequently, even a simple system such as eq(1) 
tends to behave as an extended system ( random scales behaving as 
independent degrees of freedom), when the system is allowed 
to evolve over longer and longer time scales. The scale free extension of Calculus, 
and hence of classical (and quantal ) dynamics, therefore, appears to  provide a natural 
framework for a deeper analytical treatment of a large body of natural processes 
which exhibit the phenomena of self-organised criticality \cite{sc}. We note further 
that the present formalism is also expected to offer a new approach in dealing with 
quantum-classical transitions and other related issues \cite{dec}, since the scale free 
extension can already be considered to be a non-classical extension of the classical 
dynamics.

4. It turns out that the Cantorian spacetime theory of El Naschie [8-12] 
has already anticipated some of the salient features of this extended {\em dynamical} 
Calculus. We mention, in particular, of (i) a duality (inversion $t\rightarrow 1/t) $
transformation used to obtain the Hausdorff dimension of the $\mathcal{E}^{(\infty)}$ 
space as $\phi^{-3}=4+\phi^3=4.2360679...$ \cite{eln1} from the Fisher's 
scaling relation involving low energy critical exponents (see also \cite{eln2}) 
and (ii) the non-locality in the sense that " a 'point' particle may be at two different 
spatial 'locations' at the same 'time' "~\cite{eln3}. This non-locality is a 
consequence of the indistinguishability between intersection and union in 
the underlying $\mathcal{E}^{(\infty)}$ space. It is interesting that gamma 
distributions appear naturally in the probabilistic theory of such a spacetime 
\cite{eln4}. Further, the scaling exponents $\alpha$ of 1/$f^\alpha$ noise 
fluctuations observed in semiconductor materials and quasi-crystals are related to 
the Hausdorff dimensions of  backbone Cantor sets in  ${\mathcal{E}^{n}}$ 
with mean topological dimensions $n$=4 and 5 respectively \cite{eln5}.

5. The singular role of the golden mean in the present dynamical formalism of 
Calculus suggests that a signature of the golden mean must be imprinted 
deep at the heart of every dynamical process in Nature.


\begin{thebibliography}{99}
\bibitem{pr} W H Press, Flicker noises in Astronomy and elsewhere, 
{\em Comments Astrophys.}, 1978, {\bf 7} 103-119.   
\bibitem{sc} P Bak, C Tang, W Wiesenfeld, Self-organized criticality, 
{\em Phys. Rev. A}, 1988, {\bf 38} 364-374.
\bibitem{dp} D P Datta, A new class of scale free solutions to linear ordinary 
differential equations and the universality of the Golden Mean 
${\frac{\sqrt{5}-1}{2}}=0.618033\ldots$, {\em Chaos, Solitons , Fractals} (2002), to 
appear.
\bibitem{ns} A Robinson, {\em Nonstandard analysis}, North-Holland, Amsterdam, (1966). 
\bibitem{mw} M Wolf, 1/$f$ noise in the distribution of prime numbers, 
{\em Physcia A}, 1997, {\bf 241}, 493-499.
\bibitem{pl} M Planat, 1/$f$ noise, the measurement of time and number theory, 
{\em Fluctuation and Noise Lett.}, 2001, {\bf 1}, R65-74.
\bibitem{sl} A M Selvam, Universal quantification for deterministic chaos, 
{\em Applied Math. Modelling}, 1993, {\bf 17}, 642-649.
\bibitem{eln1} M. S. El Naschie, Fisher's scaling and dualities at high energy in 
$\mathcal{E}^{(\infty)}$ spaces, {\em Chaos, Solitons \& Fractals} 2001, {\bf 12}, 
1557-1561.
 \bibitem{eln2} M. S. El Naschie, On 't Hooft dimensional regularisation in 
$\mathcal{E}^{(\infty)}$ spaces, {\em Chaos, Solitons \& Fractals}, 2001,  {\bf 12}, 
851-858. 
 \bibitem{eln3} M. S. El Naschie, On the irreducibility of spatial ambiguity in quantum 
physics,
{\em Chaos, Solitons \& Fractals}, 1998,  {\bf 6}, 913-919. 
 \bibitem{eln4} M. S. El Naschie, Remarks on superstrings, fractal gravity, Nagasawa's 
diffusion and Cantorian spacetime, {\em Chaos, Solitons \& Fractals}, 1997,  {\bf 8}, 
1873-1886. 
 \bibitem{eln5} M. S. El Naschie, Penrose tiling, semi-conduction and Cantorian 
1/$f^\alpha$ 
spectra in four and five dimensions, {\em Chaos, Solitons \& Fractals}, 1993,  {\bf 3}, 
489-491.
 \bibitem{dp1} D P Datta, Duality and scaling in quantum mechanics, 
{\em Phys. Lett. A}, 1997, {\bf  233} 274 -280.
\bibitem{pic} G F Simmons, {\em Differential equations with applications 
and historical notes}, McGraw Hill, New York, 1972. 
\bibitem{fz} G J Klir and B Yuan {\em Fuzzy Sets and Fuzzy Logic}, Prentice-Hall of 
India, 
New Delhi, 2000.
\bibitem{ot} E Ott  {\em Chaos in dynamical systems} Cambridge University
Press, Cambridge, 1993.
\bibitem{dec} D Guilini, E Joos, C Kiefer, J Kupsch, I -O Stamatescu and H D Zeh, 
{\em Decoherence and the appearance of a classical world in quantum theory}, 
Springer, Berlin, 1996.
\end{thebibliography}
\end{document}